\documentclass[reqno,
intlimits,
tbtags, %
titlepage,
12pt]{amsart}

\usepackage{microtype}
\usepackage{lmodern}
\usepackage[T1]{fontenc}
\usepackage[utf8]{inputenc}

\usepackage{times}

\usepackage{amssymb}%
\usepackage{mathrsfs}%

\usepackage[pagewise]{lineno}

\usepackage{hyperref}
\hypersetup{urlcolor = red, citecolor = blue}
\hypersetup{
  pdfauthor = Xuan Mao, 
  pdftitle = {Finite-time blow-up in chemotaxis systems with flux limitation}, 
  pdfsubject = Manuscript, 
  pdfkeywords = {Chemotaxis; flux limitation; blow-up; critical mass}
}

\usepackage[top = 1in,bottom = 1in,left = 1in,right = 1in]{geometry}
\usepackage[sort]{cite}

\usepackage{xcolor}

\theoremstyle{plain}%
\newtheorem{theorem}{Theorem}[section]
\newtheorem{corollary}[theorem]{Corollary}
\newtheorem{lemma}[theorem]{Lemma}
\newtheorem{proposition}[theorem]{Proposition}

\theoremstyle{definition}
\newtheorem{definition}{Definition}[section]

\theoremstyle{remark}
\newtheorem{remark}{Remark}[section]

\numberwithin{equation}{section}

\newcommand{\weakstar}{
\buildrel\ast\over\rightharpoonup
}

\DeclareMathOperator*{\loc}{loc}

\newcommand{\dd}{\;\mathrm{d}}
\renewcommand{\leq}{\leqslant}
\renewcommand{\geq}{\geqslant}

\newcommand{\ou}{\overline{U}}
\newcommand{\uu}{\underline{U}}

\providecommand{\e}{\ensuremath{\mathrm{e}}}
\providecommand{\ue}{\ensuremath{u^{(\varepsilon)}}}
\providecommand{\ve}{\ensuremath{v^{(\varepsilon)}}}
\providecommand{\Te}{\ensuremath{T_{\max}^\varepsilon}}
\providecommand{\Ue}{\ensuremath{U^{(\varepsilon)}}}
\providecommand{\Pd}{\Theta_{\delta}}

\DeclareSymbolFont{bbold}{U}{bbold}{m}{n}
\DeclareMathSymbol{\1}{\mathord}{bbold}{`1}

\DeclareSymbolFont{esint}{U}{esint}{m}{n}
\DeclareMathSymbol{\fintsymbol}{\mathop}{esint}{'037}
\def\fint{\fintsymbol\nolimits}

\allowdisplaybreaks[4]

\begin{document}

\title[Mass threshold in critical flux-limited chemotaxis models]{Mass threshold for global existence in chemotaxis systems with critical flux limitation}

\author[X.~Mao]{Xuan Mao}
\address[X.~Mao]{School of Mathematics\\ 
Hohai University\\ 
Nanjing 211100\\
Jiangsu, China}
\email[X.~Mao]{20250611@hhu.edu.cn}

\author[H.~Wang]{Hengling Wang}%
\address[H.~Wang]{School of Physical and Mathematical Sciences, Nanjing Tech University, Nanjing 211816, P. R. China}
\email{hlwang@njtech.edu.cn}

\author[J.~Yan]{Jianlu Yan}%
\address[J.~Yan]{School of Mathematics and Key Laboratory of MIIT, Nanjing University of Aeronautics and Astronautics, Nanjing 211106, P. R. China}
\email{ yanjl@nuaa.edu.cn}

\thanks{The first author is supported by ``the Fundamental Research Funds for the Central Universities'' (No.~B250201215).}

\thanks{The third author is supported by  Natural Science Foundation of Jiangsu Province (No. SBK2022044966).}

\subjclass[2020]{35B40; 35B33; 35B44; 35K65; 92C17}%

\keywords{Chemotaxis; flux limitation; global solvability; critical mass}

\begin{abstract}
This paper investigates the flux-limited chemotaxis system, proposed by Kohatsu and Senba~(2025),
\begin{equation*} 
  \begin{cases} 
u_t = \Delta u -\nabla\cdot(u|\nabla v|^{\alpha-2}\nabla v),\\ 
\:\:0=\Delta v + u,
  \end{cases} 
\end{equation*}  
posed in the unit ball of $\mathbb{R}^N$ for some $N\geq2$,
subject to no-flux and homogeneous Dirichlet boundary conditions.
Due to precedents, e.g., Tello (2022) and Winkler (2022), 
the exponent $\alpha = \frac{N}{N-1}$ is the threshold for finite-time blow-up under symmetry assumptions.
We further find that %
under the framework of radially symmetric solutions, 
the system with critical flux limitation admits a globally bounded weak solution if and only if initial mass is strictly less than $\omega_N \big(\frac{N^2}{N-1}\big)^{N-1}$, 
where $\omega_N$ denotes the measure of the unit sphere $\mathbb{S}^{N-1}$.
Asymptotic behaviors are also considered.
\end{abstract}

\maketitle

\section{Introduction}\label{sec: introduce}

In this paper, we consider the following initial boundary value problem for the  parabolic-elliptic chemotaxis model with flux limitation
\begin{equation}\label{sys: my ks flux limitation}
\begin{cases}
  u_t = \Delta u - \nabla \cdot(u|\nabla v|^{\alpha - 2}\nabla v), & x\in \Omega, t>0,\\
  0 = \Delta v + u, & x\in \Omega, t>0,\\
  \nabla u\cdot \nu - u|\nabla v|^{\alpha - 2}\nabla v\cdot \nu = v = 0, & x\in\partial\Omega, t>0,\\
u(\cdot, 0) = u_0(\cdot), & x\in\Omega,
\end{cases}
\end{equation}
posed in the unit ball $\Omega = B_1 := \{x\in\mathbb{R}^N: |x|<1\}$ with some $N\geq2$, 
where $\alpha = \frac{N}{N-1}$ and $\nu$ denotes the outward unit normal vector field to the boundary $\partial\Omega$.

This system is derived from the classical Keller-Segel model \cite{Keller1970}, which describes the oriented movement of cells in response to a chemical stimulus. Here, $u$ represents the density of cells and $v$ stands for the concentration of chemical signals. 
 We remark that the exponent $\alpha\neq1$ is interpreted as a generalization of flux limitation $\alpha=1$. 
For detailed biological backgrounds, refer to \cite{Bellomo2010,Bellomo2017,Bellomo2017a,Perthame2020} and the references therein.
 The goal of this paper is to investigate the global well-posedness of \eqref{sys: my ks flux limitation} with respect to critical flux-limited exponent $\alpha = \frac{N}{N-1}$.
 
 When $\alpha=2$ for $N=2$, the system \eqref{sys: my ks flux limitation} becomes free from limited flux and reduces to the Smoluchowski-Poisson equation
 \begin{equation}
  \label{sys: SP system}
   \begin{cases}
     u_t = \Delta u - \nabla\cdot(u\nabla v), & x\in\Omega, t>0, \\
     0 = \Delta v + u, & x\in\Omega, t>0.
   \end{cases}
 \end{equation}
 The system \eqref{sys: SP system} 
 was proposed by authors of \cite{Biler2006}.
 They found that %
 for any given radially symmetric initial datum with mass satisfying $\int_\Omega u_0\dd x \leq 8\pi$,
 the system \eqref{sys: SP system} admits a global solution.
 If $\int_\Omega u_0\dd x < 8\pi$, 
 the solution approaches to the steady state 
 \begin{equation}
  \label{eq: steady states of 8pi}
  \frac{8\lambda^2}{(1+|\lambda x|^2)^2}
  \quad \text{with } \lambda \text{ uniquely determined by } 
  \int_\Omega \frac{8\lambda^2}{(1+|\lambda x|^2)^2}\dd x = \int_\Omega u_0\dd x
 \end{equation}
in the space of essentially bounded functions, at least at an exponential rate,  
whereas if $\int_\Omega u_0\dd x = 8\pi$, 
it tends to $8\pi\delta_0$ in the measure space $\mathcal{M}(\overline{\Omega}) := (C(\overline{\Omega}))'$ at least algebraically.
The grow-up rate was improved by Kavallaris and Souplet~\cite{Kavallaris2009}.
Suzuki~\cite{Suzuki2013} showed that without symmetry assumptions, 
$8\pi$ is the critical mass for blow-up. 
Moreover, for a solution exploding in finite time, the set of blow-up points is finite and the solution concentrates exactly $8\pi$ mass at each blow-up point.

Considering chemotaxis models with flux limitation,
Negreanu and Tello~\cite{Negreanu2018} investigated the Neumann initial boundary value problem for the parabolic-elliptic system
\begin{equation}
  \label{sys: tello ks with flux limitation}
  \begin{cases}
    u_t = \Delta u - \chi\nabla\cdot(u |\nabla v|^{\alpha-2}\nabla v), & x\in\Omega, t > 0,\\
    0 = \Delta v - \fint_\Omega u_0\dd x + u, & x\in\Omega, t > 0,
  \end{cases}
\end{equation}
 and proved that %
   the nonnegative solution of \eqref{sys: tello ks with flux limitation} is globally bounded, provided that
 \(\alpha<\alpha_\ast\), where
  \begin{equation}
    \label{eq: alpha ast}
  \alpha_\ast : =
  \begin{cases}
    +\infty & \text{if } N=1, \\
    \frac{N}{N-1} & \text{if } N\geq2.
  \end{cases}
  \end{equation}
Subsequently, Tello~\cite{Tello2022} considered the nonnegative radial solutions of \eqref{sys: tello ks with flux limitation} and showed that for \(N>2\), \(\alpha\in(\alpha_\ast,2)\) 
and \(\chi>\chi_N\) with some \(\chi_N>0\), there exist initial data satisfying \(\fint_\Omega u_0\dd x > 6\) such that the system \eqref{sys: tello ks with flux limitation} admits a solution exploding in finite time. 
Recently, Kohatsu~\cite{Kohatsu2024} established the existence of finite-time blow-up for arbitrarily given $\chi > 0$, $\alpha > \alpha_\ast$ and initial mass. %
The threshold~\eqref{eq: alpha ast} is also the critical exponent regarding finite-time blow-up in chemotaxis models with regularized flux limitation \cite{Winkler2022,Marras2022,Mao2024a}, i.e., $|\nabla v|^{\alpha-2}$ replaced by $(1+|\nabla v|^2)^{\frac{\alpha-2}{2}}$.
For global boundedness in fully parabolic relatives,
we refer to \cite{Wang2019,Yan2020,Winkler2022a,Kohatsu2023} and references therein.
See the review~\cite{Winkler2025} for more related results.

Very recently, Kohatsu and Senba~\cite{Kohatsu2025} studied self-similar solutions to the Cauchy problem of \eqref{sys: my ks flux limitation} and particularly constructed explicitly a family of radially symmetric weak solutions $(X_\lambda, Y_\lambda)$ to the stationary problem of critical flux-limited chemotaxis model
\begin{equation}
  \label{sys: Cauchy stationary problem}
  \begin{cases}
    0 = \Delta X - \nabla\cdot (X|\nabla Y|^{\frac{1}{N-1}-1}\nabla Y), &x\in \mathbb{R}^N,\\
    0 = \Delta Y + X, & x\in\mathbb{R}^N,
  \end{cases}
\end{equation}
where 
\begin{equation}
  \label{sym: steady states of Cauchy problem}
  X_\lambda := \frac{N^{2N-1}}{(N-1)^{N-1}}\frac{\lambda^N}{(1+|\lambda x|^{\frac{N}{N-1}})^{N}}
  \quad \text{for } x\in\mathbb{R}^N \text{ and } \lambda > 0.
\end{equation}
Their work encourages us to reveal the mass threshold phenomenon concerning finite-time chemotactic collapse in critical flux-limited settings.

\subsection*{Main results and ideas}

Since the critical mass $8\pi$ of finite-time blow-up in the system~\eqref{sys: SP system} is exactly the mass of Aubin-Talenti bubble~\eqref{eq: steady states of 8pi}, we are motivated to show 
\begin{equation}
  \label{sym: critical mass}
  m_c := \omega_N \Big(\frac{N^2}{N-1}\Big)^{N-1},\quad \text{for } N = 2,3,\ldots
\end{equation}
is the mass threshold of global solvability in the critical system~\eqref{sys: my ks flux limitation}.
Indeed, the value $m_c$ is exactly the mass of $X_\lambda$ for each $\lambda > 0$.
Inductively integrating by parts, 
we have 
\begin{equation*}
  \begin{split} 
  1 &= \int_0^\infty \frac{\dd\rho}{(1+\rho)^2} 
  = \frac{\rho}{(1+\rho)^2}\biggm|_0^\infty - (-2)\int_0^\infty \frac{\rho\dd\rho}{(1+\rho)^3}\\
    &= 2\int_0^\infty \frac{\rho\dd\rho}{(1+\rho)^3}
    = 3\int_0^\infty \frac{\rho^2\dd\rho}{(1+\rho)^4} \\
    &= k\int_0^\infty\frac{\rho^{k-1}\dd\rho}{(1+\rho)^{k+1}}
    \quad \text{for all } k\in\mathbb{N}.
  \end{split}
\end{equation*}
So, for integers $N \geq 2$ and positive constants $\lambda > 0$, 
\begin{equation*}
  \begin{split}
    \int_{\mathbb{R}^N}X_\lambda\dd x
    &= \frac{Nm_c}{\omega_N}\int_{\mathbb{R}^N} \frac{\lambda^N\dd x}{(1+|\lambda x|^{\frac{N}{N-1}})^N}\\
    &= \int_0^\infty\frac{m_c\dd r}{(1+r^{\frac{1}{N-1}})^N}
    = m_c(N-1)\int_0^\infty\frac{\rho^{N-2}\dd\rho}{(1+\rho)^N} 
    = m_c,
  \end{split}
\end{equation*} 
as claimed.

Our main results, summarized below, are that under symmetry assumptions, 
the dynamics of the system \eqref{sys: my ks flux limitation} depend only on initial mass.
Here and below, let $\alpha = \frac{N}{N-1}$ and $m_c$ be given by \eqref{sym: critical mass}.

\begin{theorem}
  \label{thm: main}
  Let $N\geq2$ and $\Omega=B_1\subset\mathbb{R}^N$.
  Then for $m>0$
  whenever %
\begin{equation} 
  \label{h: initial conditions}
  u_0\in C^0(\overline{\Omega}) 
  \text{ is nonnegative and radially symmetric with } 
  \textstyle\int_\Omega u_0\dd x = m,
\end{equation}
\begin{enumerate}
\item \upshape{(Theorem~\ref{thm: finite-time blow-up})} if $m > m_c$, the corresponding solution $(u,v)$ of \eqref{sys: my ks flux limitation} blows up in finite time 
\begin{equation*} 
T_{\max} 
\leq T^\star 
:= \frac{1}{2N} \bigg(\Big(\frac{m}{m_c}\Big)^{\frac{1}{N-1}}-1\bigg)^{-1};
\end{equation*}
  \item \upshape{(Theorem~\ref{thm: global boundedness})} if $m < m_c$,
the corresponding solution $(u,v)$ of \eqref{sys: my ks flux limitation} is globally bounded and approaches to the steady state $X_\lambda$ in the topology of $L^\infty(\Omega)$,
where $\lambda > 0$ is uniquely determined by $\int_\Omega X_\lambda\dd x = m$;
\item \upshape{(Theorem~\ref{thm: infinite-time chemotactic collapse})} if $m=m_c$, the corresponding solution $(u,v)$ of \eqref{sys: my ks flux limitation} exists globally and aggregates complete mass at the center of the domain in infinite time,
\end{enumerate}
where the solution $(u,v)$ defined in $\overline{\Omega}\times[0,T_{\max})$ is given in Proposition~\ref{prop: local existence}.
\end{theorem}

Following precedents \cite{Jaeger1992,Biler2006,Kavallaris2009,Negreanu2018,Tello2022}, 
we introduce 
the accumulated density 
\begin{equation*}
  U := \int_0^{\xi^{\frac{1}{N}}}r^{N-1}u(r,t)\dd r 
  \quad \text{for } (\xi,t)\in(0,1)\times[0,T_{\max}),
\end{equation*}
and obtain our results via the study on boundedness and blow-up of the scalar parabolic equation
\begin{equation*}
  U_t = N^2\xi^{2-\frac{2}{N}}U_{\xi\xi} + N\xi^{1-\frac{2}{N}}U^{\frac{1}{N-1}}U_\xi,
\end{equation*}
subject to Dirichlet boundary conditions 
\begin{equation*}
  U(0,t) = 0
  \quad\text{and}\quad 
  U(1,t)= \frac{m}{\omega_N}
  \quad \text{for } t\in(0,T_{\max}).
\end{equation*}
Analogous to the analysis in \cite{Nagai1995,Winkler2018,Winkler2022,Wang2023a,Ahn2023}, 
we utilize the 2nd-moment function 
\begin{equation*}
  \begin{split} 
  \phi(U) &:= \frac{2\omega_N}{N}\int_0^1 \xi^{\frac{2}{N}-1} U\dd \xi
  = 2\omega_N\int_0^1\rho\int_0^\rho r^{N-1}u(r,t)\dd r\dd\rho
  \quad \text{for } t\in[0,T_{\max}),\\
  &\:= \omega_N \int_0^1\rho^{N-1}u(\rho,t)\dd\rho 
  - \omega_N\int_0^1\rho^{N+1}u(\rho,t)\dd\rho 
  = m
  - \int_\Omega |x|^2u(x,t)\dd x
  \end{split}
\end{equation*}
to study finite-time singularity formation.
The evolution of moment function is governed by  
\begin{equation*}
  \begin{split}
    \phi'(U) 
    \geq %
    2Nm\bigg(\Big(\frac{m}{m_c}\Big)^{\frac{1}{N-1}} - 1\bigg),
  \end{split}
\end{equation*}
which entails that solutions cannot be global for initial data with supercritical mass, since $\phi(U)$ is essentially a finite quantity along the solution trajectory.

Thanks to the comparison principle,
we construct a family of stationary supersolutions to show global boundedness in subcritical settings,
where these supersolutions are exactly accumulated densities of stationary solutions~\eqref{sym: steady states of Cauchy problem}. 
By contradiction, 
global existence in the case of critical mass is shown via $\varepsilon$-regularity and the strong maximum principle.
We prove both globally asymptotic stability of solutions with subcritical mass and infinite-time blow-up of solutions with critical mass by dissipation functionals.

This paper is organized as follows. 
Section~\ref{sec: local existence} is devoted to local solvability of bounded weak solutions. 
We show Theorem~\ref{thm: main} (1) in Section~\ref{sec: finite-time blow-up}.
Building upon properties of steady states discussed in Section~\ref{sec: steady states},
we show Theorem~\ref{thm: main} (2) and (3) in Section~\ref{sec: global boundedness} and Section~\ref{sec: infinite-time blow-up}, respectively.

\section{Local existence of bounded weak solutions} 
\label{sec: local existence}

Let us first introduce the natural concepts of bounded weak solutions.

\begin{definition}
  Let $T\in(0,\infty]$.
  A pair $(u,v)$ of functions: $\Omega\times(0,T)\mapsto [0,\infty)$ is called a weak solution of \eqref{sys: my ks flux limitation} in $\Omega\times(0,T)$, 
  if 
  \begin{equation*}
    \begin{split} 
    u&\in L^\infty_{\loc}([0,T);L^\infty(\Omega))\cap L^2_{\loc}([0,T); W^{1,2}(\Omega)),\\
    v&\in L^2_{\loc}([0,T); W^{1,2}(\Omega)),
    \end{split}
  \end{equation*}
  and the identities 
  \begin{equation}
    \label{eq: u weak solution}
    \begin{split} 
    \int_0^{T}\int_\Omega \nabla u\nabla \phi \dd x\dd t
    &- \int_0^{T}\int_\Omega u|\nabla v|^{\frac{1}{N-1}-1}\nabla v \nabla \phi \dd x\dd t\\
    &= \int_\Omega u_0\phi(\cdot,0)\dd x
    + \int_0^{T}\int_\Omega u\phi_t \dd x \dd t
    \end{split}
  \end{equation}
  and 
  \begin{equation*}
    \int_0^T\int_\Omega\nabla v \cdot \nabla\psi \dd x\dd t  
    = \int_0^T\int_\Omega u\psi\dd x \dd t 
  \end{equation*}
  hold for all  $\phi,\psi\in C^\infty_0(\overline\Omega\times[0,T))$.  
\end{definition}

The local existence and uniqueness of bounded weak solutions to \eqref{sys: my ks flux limitation}  can be asserted by adopting well-established fixed-point arguments (cf.~\cite[Proposition~2.1]{Fuhrmann2022} and \cite[Theorem~2.1]{Li2023}) for parabolic-elliptic chemotaxis systems under no-flux boundary conditions for cell density $u$.  

\begin{proposition}
\label{prop: local existence}
  Let $N\geq2$ and $\Omega\subset\mathbb{R}^N$ be the unit ball, 
  and assume that \(u_0\in C^0(\overline{\Omega})\) is nonnegative and radially symmetric. 
  Then there exist $T_{\max}\in(0,\infty]$ and a uniquely determined pair $(u,v)$ of nonnegative functions 
  \begin{equation*}
    \begin{cases}
      u\in C^0(\overline{\Omega}\times[0,T_{\max}))\cap L^2_{\loc}([0,T_{\max}); W^{1,2}(\Omega)), \\
      v\in \cap_{q>N}C^0([0,T_{\max});W^{2,q}(\Omega)\cap W^{1,q}_0(\Omega)),
    \end{cases}
  \end{equation*} 
  such that $(u,v)$ is a weak solution to \eqref{sys: my ks flux limitation} in $\Omega\times(0, T_{\max})$, 
  that
  \begin{equation}\label{eq: mass conservation}
    \int_\Omega u(\cdot, t) \dd x = \int_\Omega u_0 \dd x \quad\text{for all } t\in(0,T_{\max}),
  \end{equation}
  that
  \begin{equation}\label{eq: extensibility criterion}
    \text{if } T_{\max}<\infty, \text{ then } \limsup_{t\nearrow T_{\max}}\|u(\cdot, t)\|_{L^\infty(\Omega)} = \infty,
  \end{equation}
  and that $u(\cdot, t)$ and $v(\cdot, t)$ are radially symmetric for each $t\in(0, T_{\max})$.
  Furthermore, 
  \begin{align*}
    u &\in C^{\vartheta,\vartheta/2}(\overline{\Omega}\times(0,T_{\max}))
    \cap 
    C^{2,1}(\overline{\Omega}\setminus\{0\}\times(0,T_{\max})) \\
    v &\in C^{2,0}(\overline{\Omega}\times(0, T_{\max}))
  \end{align*}
  is positive for all $(x,t)\in\overline{\Omega}\times(0,T_{\max})$, which solves \eqref{sys: my ks flux limitation} classically in $\overline{\Omega}\setminus\{0\}\times[0,T_{\max})$.
\end{proposition}

\begin{proof}
  The existence, uniqueness and extensibility criterion of bounded weak solutions to \eqref{sys: my ks flux limitation} can be inferred from fixed-point theories as the lines in \cite[Proposition~2.1]{Fuhrmann2022}.
  To show mass conservation~\eqref{eq: mass conservation}, 
  we introduce $\Upsilon\in C^\infty(\mathbb{R})$, a nonincreasing function with properties
  $\Upsilon\equiv1$ for $r\leq0$ and $\Upsilon\equiv0$ for $r\geq1$.
  Then for $\tau\in(0,T_{\max})$ substituting $\phi(x,t) := \Upsilon(\delta^{-1} (t-\tau))$ into \eqref{eq: u weak solution},
  we get 
  \begin{equation*} 
    0 = \int_\Omega u_0\dd x 
    + \fint_\tau^{\tau+\delta} \Upsilon'(\delta^{-1}(t-\tau))\int_\Omega u(x,t)\dd x \dd t
    \quad \text{for all } \delta\in(0,T_{\max}-\tau),
  \end{equation*}
  which implies 
  \begin{equation*}
    \bigg|\int_\Omega (u_0 - u(x,\tau))\dd x\bigg|
    \leq \sup_{t\in(\tau,\tau+\delta)}\bigg|\int_\Omega(u(x,t)-u(x,\tau))\dd x\bigg|
    \quad\text{for all } \delta\in(0,T_{\max}-\tau).
  \end{equation*}
  We have $\int_\Omega (u_0 - u(x,\tau))\dd x = 0$ by $\delta\downarrow0$ due to continuity of $u$.
  Radial symmetry is consequences of uniqueness and the exact form of \eqref{sys: my ks flux limitation}.
  The positivity of $u$ follows from Harnack inequality, and H\"older continuity can be derived from \cite[Theorem~1.3]{Porzio1993}. %
  Thus, by standard elliptic regularity theories,
  \begin{equation*}
    v\in C^{2+\vartheta,1+\vartheta/2}(\overline{\Omega}\times(0,T_{\max})) 
  \quad\text{for some } \vartheta\in(0,1). 
  \end{equation*}
By virtue of symmetry, we may obtain $\nabla v$ by integrating the second equation $0 = (v_rr^{N-1})_r + ur^{N-1}$ over $(0,\rho)\subset(0,1)$ accordingly  
  \begin{equation}
    \label{eq: grad v =}
    - v_r(\rho,t) = \rho^{1-N}\int_0^\rho ur^{N-1}\dd r > 0 
    \quad \text{for all } (\rho,t)\in(0,1]\times(0,T_{\max})
  \end{equation} 
  and
  \begin{equation*}
    |\nabla v|^{\frac{N}{N-1} - 2} \in C_{\loc}^{1+\vartheta, 1+\vartheta/2}(\overline{\Omega}\setminus\{0\}\times(0,T_{\max})).
  \end{equation*} 
  Thanks to \cite{Lieberman1987}, H\"older gradient estimates of $u$ in $\overline{\Omega}\setminus\{0\}\times(0,T_{\max})$ can be obtained.
  Therefore, we may employ the standard Schauder theories to show that
  \begin{equation*}
    u \in C^{2,1}(\overline{\Omega}\setminus\{0\}\times(0,T_{\max}))
  \end{equation*}
  which together with $v$ solves \eqref{sys: my ks flux limitation} in $\overline{\Omega}\setminus\{0\}\times(0,T_{\max})$ classically.
\end{proof}

Following the precedents in~\cite{Jaeger1992,Tao2017,Tao2025}, 
we shall analyze a scalar parabolic problem satisfied by accumulated densities to detect blow-up and establish boundedness under symmetry assumptions.

\begin{lemma} \label{le: mass distribution function}
  Let the pair $(u,v)$ of radially symmetric functions be a bounded weak solution to the system~\eqref{sys: my ks flux limitation} defined in $\overline\Omega\times[0,T_{\max})$ given by Proposition~\ref{prop: local existence}.
  Then, the mass distribution function
\begin{equation}
  \label{sym: U mass distribution function}
  U(\xi, t) := \int_0^{\xi^\frac1N}u(r,t)r^{N-1}\dd r
  \quad \text{for } (\xi,t)\in[0,1]\times[0,T_{\max}),
\end{equation}
satisfying 
\begin{equation} 
  \label{eq: U regularity}
  U\in C^0([0,T_{\max}); C^1([0,1]))\cap C^{2,1}((0,1]\times(0,T_{\max})) 
\end{equation}
and 
\begin{equation}
  \label{eq: U > 0}
  U > 0\quad\text{and}\quad U'>0 \quad \text{for all }(\xi,t)\in(0,1]\times(0,T_{\max}),
\end{equation}
solves the scalar parabolic problem 
\begin{equation}
  \label{sys: partial mass with flux limitation}
\begin{cases}
  \mathcal{P}(U) = 0 & (\xi,t)\in(0,1)\times(0,T_{\max}), \\
  U(0, t)=0, U(1, t) = \frac{m}{\omega_N}, & t\in(0,T_{\max}), \\
  U(\cdot, 0) = U_0(\cdot) & \xi\in(0,1),
\end{cases}
\end{equation}
classically,
where
\begin{equation}
\label{sym: mathcal P}
\mathcal{P}(U) := U_t - N^2\xi^{2-\frac2N}U_{\xi\xi} 
- N\xi^{1-\frac2N}U^{\frac{1}{N-1}}U_\xi
\end{equation}
and
\begin{equation*}
  \label{eq: U0 ks flux limitation}
 U_0(\xi) 
 := \int_0^{\xi^{\frac1N}}u_0(r)r^{N-1}\dd r 
 \quad\text{for } \xi\in[0,1].
\end{equation*}
\end{lemma}

\begin{proof}
  Using  
  \begin{equation*}
    U_\xi = \frac{u(\xi^{\frac{1}{N}},t)}{N},
    \quad 
    U_{\xi\xi} = \frac{u_r(\xi^{\frac{1}{N}},t)}{N^2}\xi^{\frac{1}{N}-1}, 
  \end{equation*}
  and \eqref{eq: grad v =}, i.e., $-v_r(\xi^{\frac{1}{N}},t) = \xi^{\frac{1}{N}-1}U$,
  we calculate %
  \begin{equation}
    \label{eq: d/dt outer region}
    \begin{split} 
    \frac{\dd }{\dd t}\int^1_{\xi^{\frac{1}{N}}}u(r,t)r^{N-1}\dd r
    &=\int_{\xi^{\frac{1}{N}}}^1(u_rr^{N-1}-u|v_r|^{\frac{1}{N-1}-1}v_rr^{N-1})_r\dd r\\
    &= - N^2\xi^{2-\frac{2}{N}}U_{\xi\xi} 
    - N\xi^{1-\frac{2}{N}}U^{\frac{1}{N-1}}U_\xi
    \end{split}
  \end{equation}
  for $(\xi,t)\in(0,1)\times(0,T_{\max})$, 
  where we have used the no-flux boundary conditions.
Thanks to mass identity~\eqref{eq: mass conservation}, 
i.e., $\int_0^1ur^{N-1}\dd r = U(1,t) \equiv \frac{m}{\omega_N}$ for all $t\in(0,T_{\max})$ and  
\begin{equation*}
  \frac{\dd}{\dd t}\int_0^1ur^{N-1}\dd r \equiv 0
  \quad\text{for all } t\in(0,T_{\max}), 
\end{equation*}
\eqref{eq: d/dt outer region} is reduced to $\mathcal{P}(U) = 0$ for all $(\xi,t)\in(0,1)\times(0,T_{\max})$. 
\end{proof}

\section{Finite-time blow-up}
\label{sec: finite-time blow-up}

This section is devoted to show that the solution explodes in finite time, 
for arbitrary radially symmetric initial datum with supercritical mass.

\begin{theorem}
  \label{thm: finite-time blow-up}
  Let $N\geq2$ and $\Omega = B_1\subset\mathbb{R}^N$.
  If $u_0\in C^0(\overline{\Omega})$ complying with \eqref{h: initial conditions} has subcritical mass, i.e., 
  \begin{equation}
    \label{h: supercritical mass}
    m = \int_\Omega u_0\dd x > m_c,
  \end{equation}
  then the corresponding solution given by Proposition~\ref{prop: local existence} blows up in finite time $T_{\max} \leq T^\star$,
  where $T^\star < \infty$ is defined in Theorem~\ref{thm: main}.
\end{theorem}

Our main tool is a moment inequality.

\begin{lemma}
  \label{le: moment inequality}
  Let $U$ be the accumulated density given by Lemma~\ref{le: mass distribution function}.
  Define  
  \begin{equation}
    \label{sym: moment}
    \psi(t) := \int_0^1U\xi^{\frac{2}{N}-1}\dd\xi  
    \quad \text{for } t\in[0,T_{\max}).
  \end{equation}
  Then 
  \begin{equation}
    \label{eq: moment inequality}
    \psi(t) - \psi(0)
    \geq \frac{N^2m}{\omega_N} 
    \bigg(\Big(\frac{m}{m_c}\Big)^{\frac{1}{N-1}}-1\bigg) t 
    \quad \text{for all } t\in(0,T_{\max}).
  \end{equation}
\end{lemma}

\begin{proof}
  Thanks to \eqref{eq: U regularity}, \eqref{eq: U > 0} and $U(0,t) = 0$ for all $t\in[0,T_{\max})$,
  we compute 
  \begin{equation}
    \label{eq: odi of psi}
    \begin{split}
      \psi'
      &= \int_0^1 N^2\xi U_{\xi\xi} \dd\xi 
      + \int_0^1 NU^{\frac{1}{N-1}}U_\xi\dd\xi\\
      &= \int_0^1\Big(N^2\xi U_\xi - N^2U + (N-1)U^{\frac{N}{N-1}}\Big)_\xi\dd\xi\\
      &= N^2U_\xi(1,t) 
      - \frac{N^2m}{\omega_N} 
      + (N-1)\bigg(\frac{m}{\omega_N}\bigg)^{\frac{N}{N-1}}\\
      &\geq \frac{N^2m}{\omega_N} 
      \bigg(\Big(\frac{m}{m_c}\Big)^{\frac{1}{N-1}}-1\bigg)
    \end{split}
  \end{equation}
  for all $t\in(0,T_{\max})$.
Integrating \eqref{eq: odi of psi} 
we obtain \eqref{eq: moment inequality}.
\end{proof}

We are in a position to show Theorem~\ref{thm: finite-time blow-up}.

\begin{proof}
  [Proof of \upshape{Theorem~\ref{thm: finite-time blow-up}}]
  Since 
  \begin{equation*}
    \psi \in \bigg(0, \frac{Nm}{2\omega_N}\bigg)
    \quad \text{for all } t\in[0,T_{\max}),
  \end{equation*}
  according to the definition~\eqref{sym: moment},
  Lemma~\ref{le: moment inequality} warrants that  
  \begin{equation*}
    \frac{N^2m}{\omega_N} 
      \bigg(\Big(\frac{m}{m_c}\Big)^{\frac{1}{N-1}}-1\bigg)t
    < \frac{Nm}{2\omega_N} 
    \quad \text{for all } t\in(0,T_{\max}).
  \end{equation*}
  It follows from \eqref{h: supercritical mass} that 
  \begin{equation*}
    T_{\max} \leq \frac{1}{2N} \bigg(\Big(\frac{m}{m_c}\Big)^{\frac{1}{N-1}}-1\bigg)^{-1},
  \end{equation*}
  and completes the proof.
\end{proof}

\begin{remark}
  The generalized moment method gives a simple proof for generic finite-time blow-up. 
  We would like to mention that based on maximum principle (see Lemma~\ref{le: comparison principle} below), 
  gradient exploding subsolutions $\psi_\lambda$ can be constructed to show existence of finite-time blow-up for a restricted class of initial data with supercritical mass,
  where %
  \begin{equation*}
    \psi_\lambda(s,t) := \e^{\tau} W_\lambda(\tau^{-3(N-1)}s) 
    \quad \text{for all } (s,t)\in[0,\infty)\times[0,\varepsilon/\delta) \text{ and } \lambda > 0,
  \end{equation*}
  with $\tau := \varepsilon - \delta t$, 
  \begin{equation*}
    \varepsilon := \ln\frac{m+\omega_NA}{2\omega_NA}
    \quad\text{and}\quad 
    \delta := \frac{\lambda^{-1}\gamma^{-2}}{N-1}.
  \end{equation*}
  Such construction is motivated by \cite{Mao2025}, where the case $N=2$ has been shown. 
\end{remark}

\section{Steady states}
\label{sec: steady states}

In contrast to the supercritical mass case, 
we shall show that the system \eqref{sys: my ks flux limitation} with critical flux limitation admits stationary solutions, 
which have subcritical mass.

\begin{proposition}
  \label{prop: steady states}
  Let $N\geq2$.
  Assume that the nonnegative function $W\in C^0([0,\infty))\cap C^2((0,\infty))$ is a classical solution to the one-point boundary value problem 
    \begin{equation}
      \label{sys: one-point boundary value problem}
      \begin{cases}
        N\xi W_{\xi\xi} + W^{\frac{1}{N-1}}W_\xi = 0, & \xi > 0,\\
        W(0) = 0, 
      \end{cases}  
    \end{equation}
  such that 
  \begin{equation}
    \label{eq: W not euqiv 0}
    W\not\equiv0.
  \end{equation}      
  Then there exists $\lambda > 0$ such that 
  \begin{equation*}
    W = W_\lambda 
    \quad \text{for } \xi\in[0,\infty), 
  \end{equation*}
  where $W_\lambda := W_0(\lambda^N\xi)$ and 
  \begin{equation*}
    W_0 := \Big(\frac{N^2}{N-1}\Big)^{N-1}\xi\big(1+\xi^{\frac{1}{N-1}}\big)^{1-N}
    \quad \text{for } \xi \geq 0.
  \end{equation*}
  For each $\lambda > 0$, $W_\lambda$ complying with  
  \begin{equation*}
    W_\lambda \to A := \Big(\frac{N^2}{N-1}\Big)^{N-1} \quad\text{as } \xi\to\infty
  \end{equation*}
  is strictly increasing and concave.
\end{proposition}
 
\begin{proof}
  Existence. 
We compute 
\begin{equation}
  \label{eq: W lambda xi}
  \begin{split} 
  W_{\lambda\xi} &= \lambda^N A\Big(1+(\lambda^N\xi)^{\frac{1}{N-1}}\Big)^{1-N}
  - \lambda^N A(\lambda^N\xi)^{\frac{1}{N-1}}\Big(1+(\lambda^N\xi)^{\frac{1}{N-1}}\Big)^{-N}\\
  & = \lambda^N A\Big(1+(\lambda^N\xi)^{\frac{1}{N-1}}\Big)^{-N}
  \end{split}
\end{equation}
and 
\begin{equation*}
  \begin{split} 
  N\xi W_{\lambda\xi\xi} 
  &= - \frac{\lambda^N AN^2}{N-1}\Big(1+(\lambda^N\xi)^{\frac{1}{N-1}}\Big)^{-N-1}(\lambda^N\xi)^{\frac{1}{N-1}} \\
  &= - \frac{N^2}{N-1} \Big(1+(\lambda^N\xi)^{\frac{1}{N-1}}\Big)^{-1}(\lambda^N\xi)^{\frac{1}{N-1}}W_{\lambda\xi} 
  = - W_\lambda^{\frac{1}{N-1}}W_{\lambda\xi}
  \end{split}
\end{equation*}
for all $\xi > 0$, which implies $W_\lambda$ is indeed a solution to \eqref{sys: one-point boundary value problem} for each $\lambda > 0$.

Uniqueness up to dilations.
Using Newton-Leibniz formula, we compute 
\begin{equation*}
  \begin{split}
    W(\xi) - W(\eta) 
    &= \int_\eta^\xi W_\xi \dd \zeta 
    = W_\xi(\xi)\xi 
    - W_\xi(\eta)\eta 
    - \int_\eta^\xi \zeta W_{\xi\xi}\dd \zeta \\
    &= W_\xi(\xi)\xi
    - W_\xi(\eta)\eta 
    + \frac{1}{N}\int_\eta^\xi W^{\frac{1}{N-1}}W_\xi\dd \zeta \\
    &= W_\xi(\xi)\xi
    - W_\xi(\eta)\eta 
    + \frac{N-1}{N^2}W^{\frac{N}{N-1}}(\xi) 
    - \frac{N-1}{N^2}W^{\frac{N}{N-1}}(\eta)
  \end{split}
\end{equation*}
for $0 < \eta < \xi < \infty$. 
Taking $\eta\downarrow0$, we obtain the limit $K := \lim_{\eta\downarrow0}W_s(\eta)\eta$ exists and 
\begin{equation*}
  K = W_\xi(\xi)\xi
  + \frac{N-1}{N^2}W^{\frac{N}{N-1}}(\xi) 
    - W(\xi)
    \quad \text{for all } \xi \in(0, \infty),
\end{equation*} 
where we have used the boundary condition $W(0) = 0$.
Since $W$ cannot be a continuous function up to $0$ in the case of $K \neq 0$, which is a contradiction, 
we have $K = 0$ and 
thus,
\begin{equation}
  \label{sys: ode of W}
  \xi W_\xi = W - \frac{N-1}{N^2}W^{\frac{N}{N-1}}
  = W\bigg(1 - \frac{N-1}{N^2}W^{\frac{1}{N-1}}\bigg)
  \quad \text{for } \xi \in (0, \infty).
\end{equation}
Write
\begin{equation*}
  f := W^{\frac{1}{N-1}}
\end{equation*}
for brevity.
Then \eqref{sys: ode of W} is reduced to 
\begin{equation}
  \label{eq: equation of f}
  (N-1)\xi f_\xi 
  = f\bigg(1 - \frac{N-1}{N^2}f\bigg)
  \quad \text{for } \xi \in (0, \infty), \text{ with } f(0) = 0.
\end{equation}
If there exists $\zeta\in(0,\infty)$ such that $f(\zeta) = 0$,
then by uniqueness theorem of ordinary differential equation, we have $f\equiv 0$,
which is incompatible with our assumption \eqref{eq: W not euqiv 0}.
So $f > 0$ for all $\xi\in(0,\infty)$ and we may define the connected set 
\begin{equation*}
  \mathcal{S} := \bigg\{\xi\in(0,\infty): f(\zeta)\in\bigg(0,\frac{N^2}{N-1}\bigg)\text{ for all }\zeta\in(0,\xi]\bigg\},
\end{equation*} 
which is nonempty. %
Rearranging \eqref{eq: equation of f}, we get 
\begin{equation*}
  \frac{f_\xi}{f} + \frac{\frac{N-1}{N^2}f_\xi}{1-\frac{N-1}{N^2}f} = \frac{1}{(N-1)\xi} 
  \quad \text{for } \xi \in \mathcal{S}.
\end{equation*}
Integrating the equation above over $(0,\xi)\subset\mathcal{S}$ yields 
\begin{equation*}
  \frac{f}{1-\frac{N-1}{N^2}f} = \frac{N^2}{N-1}(\lambda^N\xi)^{\frac{1}{N-1}}
  \quad \text{for } \xi \in \mathcal{S} \text{ with some } \lambda > 0.
\end{equation*}
that is, 
\begin{equation*}
  f = \frac{N^2}{N-1}\frac{(\lambda^N\xi)^{\frac{1}{N-1}}}{1+(\lambda^N\xi)^{\frac{1}{N-1}}}
  \quad \text{for } \xi \in \mathcal{S},
\end{equation*}
which implies $\mathcal{S} = (0,\infty)$ and $W$ coincides with $W_\lambda$ for some $\lambda>0$.
\end{proof}

\begin{remark}
  The stationary solutions to \eqref{sys: one-point boundary value problem} are accumulated densities of radially symmetric stationary solutions \eqref{sym: steady states of Cauchy problem} to the Cauchy problem \eqref{sys: Cauchy stationary problem} for the stationary system of critical flux-limited chemotaxis models~\cite[Theorem~1.3]{Kohatsu2025}, which can be seen from \eqref{eq: W lambda xi}.
\end{remark}

\begin{corollary}
  \label{coro: steady states}
  Let $N\geq2$.
  Assume that the nonnegative function $\phi_\ell\in C^0([0,1])\cap C^2((0,1])$ is a classical solution to the Dirichlet boundary value problem 
    \begin{equation*}
      \label{sys: Dirichlet boundary value problem}
      \begin{cases}
        N\xi \phi_{\xi\xi} + \phi^{\frac{1}{N-1}}\phi_\xi = 0, & \xi \in (0,1),\\
        \phi(0) = 0,\quad \phi(1) = \ell, 
      \end{cases}  
    \end{equation*}
for some $\ell > 0$.
Then 
  \begin{equation*}
    \ell < A 
  \end{equation*}
  and 
  \begin{equation*}
    \phi_\ell = W_\lambda,
  \end{equation*}
  where $\lambda > 0$ is uniquely determined by  
  \begin{equation}
    \label{eq: lambda determined by ell}
    W_0(\lambda^N) = \ell. 
  \end{equation}
  For each $\xi \in (0,1]$, the mapping $\ell\mapsto \phi_\ell(\xi)$ complying with 
  \begin{equation}
    \label{eq: phi ell uparrow A}
    \phi_\ell \uparrow A \quad\text{as } \ell\uparrow A,
  \end{equation}
  is strictly increasing with respect to $\ell\in(0,A)$.
\end{corollary}

\section{Global boundedness and asymptotic stability}
\label{sec: global boundedness}

This section is devoted to global boundedness and asymptotic stability of solutions with subcritical mass.

\subsection{\texorpdfstring{$\varepsilon$}{ε}-regularity}

For generality, we established an $\varepsilon$ regularity,
which will be useful for global solvability of solutions with critical mass discussed in Section~\ref{sec: infinite-time blow-up}.

\begin{proposition}
  \label{prop: varepsilon regualrity}
  Let $N\geq2$ and $\Omega = B_1\subset\mathbb{R}^N$.
  If there exists $\rho\in(0,1]$ such that 
  \begin{equation}
    \label{eq: U(rho,t) < A}
    \sup_{t\in(0,T_{\max})} U(\rho,t) < A,
  \end{equation}
  then there exists $C>0$ such that
  \begin{equation}
    \label{eq: U xi has upper bounds}
    U_\xi(\xi,t) \leq C\quad\text{for all }(\xi,t)\in(0,1)\times(0,T_{\max}).
  \end{equation}
  Particularly, for any radially symmetric initial datum $u_0\in C^0(\overline{\Omega})$ with subcritical mass $\int_\Omega u_0\dd x < m_c$, 
  the corresponding solution exists globally and remains bounded.
\end{proposition}

The proof of the above proposition consists of several lemmas.
We begin with the following comparison principle.

\begin{lemma}
  \label{le: comparison principle}
Let $N\geq2$ and $\Omega=B_1\subset\mathbb{R}^N$. 
Suppose that nonnegative and nondecreasing functions $\ou$ and $\uu$ from $C^0([0,T);C^1([0,1]))\cap C^{2,1}((0,1]\times(0,T))$ satisfy 
\begin{equation}
    \label{h: uu leq ou initial}
    \uu(\cdot,0) \leq \ou(\cdot,0)  
    \quad\text{for all } \xi\in(0,1),
\end{equation}
and for some $\Lambda>0$,
\begin{equation}
  \label{h: uu leq ou boundary conditions}
  \uu(0,t)\leq\ou(0,t)
  \quad\text{and}\quad
  \uu(1,t)\leq \ou(1,t) \leq \Lambda
  \quad \text{for all } t\in (0,T),
\end{equation}
as well as 
\begin{equation}
  \label{h: puu leq pou}
  \mathcal{P}(\uu) \leq \mathcal P(\ou) 
  \quad \text{for all } (\xi,t)\in(0,1)\times(0,T).
\end{equation}
Then 
\begin{equation*}
     \uu \leq \ou 
     \quad\text{for all }\xi\in(0,1)\text{ and }t\in(0,T).
\end{equation*}
\end{lemma}

\begin{proof}
  The proof is a minor modification of \cite[Proposition~3.1]{Kavallaris2009},
  where the case $N = 2$ has been proven.
  Let $Z := \uu - \ou$, 
  then 
  \begin{equation}
    \label{eq: Z regualrity}
    Z\in C^0([0,T);C^1([0,1]))\cap C^{2,1}((0,1]\times(0,T)).
  \end{equation} 
  By \eqref{h: puu leq pou}, we have 
  \begin{equation}
    \label{eq: Zt leq etc}
    \begin{split}
      Z_t &\leq N^2\xi^{2-\frac{2}{N}}Z_{\xi\xi} 
      + N\xi^{1-\frac{2}{N}}\uu^{\frac{1}{N-1}}\uu_\xi 
      - N\xi^{1-\frac{2}{N}}\ou^{\frac{1}{N-1}}\ou_\xi\\
      &= \xi^{1-\frac{2}{N}}\Big(N^2\xi Z_\xi 
      - N^2Z 
      + (N-1)\uu^{\frac{N}{N-1}}
      - (N-1)\ou^{\frac{N}{N-1}}\Big)_\xi 
    \end{split}
  \end{equation}
  for all $(\xi,t)\in(0,1)\times(0,T)$.
  For $\delta > 0$, we define the following $C^1$ (and piecewise $C^2$) convex approximations of the mapping $s\mapsto s_+ := \max\{s,0\}$:
  \begin{equation}
    \label{sym: Pd}
    \Pd(s) := 
    \begin{cases}
      0 & \text{if } -\infty < s \leq \delta,\\
      (2\delta)^{-1}(s-\delta)^2 & \text{if } \delta \leq s\leq 2\delta,\\
      s - 3\delta/2 &\text{if } 2\delta < s < \infty.
    \end{cases}
  \end{equation} 
  Fix $0 < \tau < \eta < T$ and $\delta \in(0,1)$.
  Then for any $\varepsilon\in(0,1)$, 
  multiplying~\eqref{eq: Zt leq etc} by $\xi^{\frac{2}{N}-1}\Pd'(Z)$, 
  integrating by parts,
  we obtain 
  \begin{equation}
    \label{eq: phi delta eta - tau}
    \begin{split} 
    &\quad \int_\varepsilon^1\xi^{\frac{2}{N}-1}\Pd(Z(\xi,\eta))\dd \xi 
    - \int_\varepsilon^1\xi^{\frac{2}{N}-1}\Pd(Z(\xi,\tau))\dd \xi 
    = \int_{\tau}^\eta \int_\varepsilon^1 \xi^{\frac{2}{N}-1}\Pd'(Z)Z_t\dd\xi\dd t \\
    &\leq \int_\tau^\eta\int_\varepsilon^1 \Pd'(Z)
    \Big(N^2\xi Z_\xi 
      - N^2Z 
      + (N-1)\uu^{\frac{N}{N-1}}
      - (N-1)\ou^{\frac{N}{N-1}}\Big)_\xi \dd\xi\dd t\\
    &= \int_\tau^\eta\bigg[\Pd'(Z)
    \Big(N^2\xi Z_\xi 
      - N^2Z 
      + (N-1)\uu^{\frac{N}{N-1}}
      - (N-1)\ou^{\frac{N}{N-1}}\Big)\bigg]^1_\varepsilon\dd t\\
    &\quad - \int_\tau^\eta\int_\varepsilon^1 \Pd''(Z)
    \Big(N^2\xi Z_\xi 
      - N^2Z 
      + (N-1)\uu^{\frac{N}{N-1}}
      - (N-1)\ou^{\frac{N}{N-1}}\Big)Z_\xi \dd\xi\dd t\\
      &=: J_1 - J_2.
    \end{split}
  \end{equation}
  Using the second inequality in \eqref{h: uu leq ou boundary conditions}, $0\leq\Pd'\leq1$, and the assumption that $\uu$ and $\ou$ are nonnegative and nondecreasing,
  we estimate 
  \begin{equation*}
    \begin{split}
      J_1 
      &\leq N^2\varepsilon\int_\tau^\eta|Z_\xi(\varepsilon,t)|\dd t 
      + \int_\tau^\eta \Pd'(Z(\varepsilon,t))\bigg(N^2Z(\varepsilon,t) + (N-1)\ou^\frac{N}{N-1}(\varepsilon,t)\bigg)\dd t \\
      &\leq N^2\varepsilon\int_\tau^\eta|Z_\xi(\varepsilon,t)|\dd t 
      + \Big(N^2\Lambda 
      + (N-1)\Lambda^{\frac{N}{N-1}} \Big)\int_\tau^\eta \Pd'(Z(\varepsilon,t))\dd t.
    \end{split}
  \end{equation*}
  In light of $\Pd'' \geq 0$ and mean value theorem,  
\begin{equation}
  \label{eq: mean value theorem}
  \begin{split}
    \bigg\lvert \xi^{\frac{N}{N-1}} - \zeta^{\frac{N}{N-1}}\bigg\rvert 
    &\leq \frac{N}{N-1}\Lambda^{\frac{1}{N-1}}|\xi-\zeta|
    \quad \text{for all } \xi,\zeta\in[0,\Lambda],
  \end{split}
\end{equation}
we estimate 
\begin{align*}
  - J_2 &\leq - \int_\tau^\eta\int_\varepsilon^1 \Pd''(Z)
  \Big( 
    - N^2Z 
    + (N-1)\uu^{\frac{N}{N-1}}
    - (N-1)\ou^{\frac{N}{N-1}}\Big)Z_\xi \dd\xi\dd t\\
    &\leq  \Big(N^2 
    + N\Lambda^{\frac{1}{N-1}}\Big)\int_\tau^\eta\int_\varepsilon^1|\Pd''(Z)Z||Z_\xi|\dd\xi\dd t.
\end{align*}
Write 
\begin{equation*}
  C:= N^2\Lambda 
  + (N-1)\Lambda^{\frac{N}{N-1}} 
  + N^2 
  + N\Lambda^{\frac{1}{N-1}}.
\end{equation*}
Then 
\eqref{eq: phi delta eta - tau} is reduced to 
\begin{align*}
    &\quad \int_\varepsilon^1\xi^{\frac{2}{N}-1}\Pd(Z(\xi,\eta))\dd \xi 
    - \int_\varepsilon^1\xi^{\frac{2}{N}-1}\Pd(Z(\xi,\tau))\dd \xi \\
    &\leq N^2\varepsilon\int_\tau^\eta|Z_\xi(\varepsilon,t)|\dd t 
    + C\int_\tau^\eta \Pd'(Z(\varepsilon,t))\dd t
    + C\int_\tau^\eta\int_\varepsilon^1|\Pd''(Z)Z||Z_\xi|\dd\xi\dd t
\end{align*}
which implies 
\begin{equation*}
  \int_0^1\xi^{\frac{2}{N}-1}\Pd(Z(\xi,\eta))\dd \xi 
    - \int_0^1\xi^{\frac{2}{N}-1}\Pd(Z(\xi,\tau))\dd \xi 
    \leq C\int_\tau^\eta\int_0^1|\Pd''(Z)Z||Z_\xi|\dd\xi\dd t
\end{equation*}
by dominated convergence theorem used in the limit $\varepsilon\downarrow0$,
thanks to \eqref{h: uu leq ou boundary conditions}, 
\eqref{eq: Z regualrity} and the definition of $Z$.
Observe that $\lim_{\delta\downarrow0}\Pd(s) = s_+$ and 
$s\Pd''(s)$ is uniform-in-$\delta$ bounded and converges to $0$ for a.e. $s\in\mathbb{R}$, as $\delta\downarrow0$.
Using $0\leq \Pd(s)\leq s_+$,
we may pass to the limit $\delta\downarrow0$ by dominated convergence theorem in the preceding inequality, and we obtain 
\begin{equation*}
  \int_0^1\xi^{\frac{2}{N}-1}Z_+(\xi,\eta)\dd \xi 
    - \int_0^1\xi^{\frac{2}{N}-1}Z_+(\xi,\tau)\dd \xi 
    \leq 0.
\end{equation*}
Letting $\tau\downarrow0$ and using \eqref{h: uu leq ou initial}, 
we conclude that $\int_0^1Z_+(\xi,\eta) \leq 0$ for all $\eta\in(0,T)$, 
and finish the proof.
\end{proof}

We need the following boundedness principle to show the global boundedness of solutions. 
Analogous properties \cite{Winkler2019,Fuhrmann2022,Mao2024} have been observed in contexts of classical solutions to parabolic-elliptic chemotaxis models without flux limitation, obtained via either Bernstein-type argument or Nagai method. 
Here, for completeness, 
we give a proof within the framework of weak solutions.
In particular, the case of $N=2$ has been considered in \cite{Kavallaris2009} via boundedness arguments on a reaction-diffusion-drift equation satisfied by the average accumulated mass function $r^{-N}\int_0^ru(\rho,t)\rho^{N-1}\dd\rho$.

\begin{lemma}
  \label{le: varepsilon regularity}
Let $N\geq2$ and $\Omega=B_1\subset\mathbb{R}^N$. 
Suppose that $u_0\in C^0(\overline{\Omega})$ satisfies \eqref{h: initial conditions} and 
\begin{equation}\label{h: U/s has upper bounds}
  \sup_{(\xi,t)\in(0,1)\times(0,T_{\max})} \frac{U(\xi,t)}{\xi} < \infty,
\end{equation}
where $T_{\max}$ is given as in Proposition~\ref{prop: local existence} and $U$ is defined in \eqref{sym: U mass distribution function}. 
Then \eqref{eq: U xi has upper bounds} holds.
\end{lemma}

\begin{proof}
  Inspired by \cite{Yan2020}, 
  we consider classical solutions $(\ue,\ve)$ of the regularized systems for $\varepsilon\in(0,1)$,
\begin{equation}
  \label{sys: regularized system}
  \begin{cases}
    \ue_{t} = \Delta \ue - \nabla\cdot(\ue(\varepsilon + |\nabla \ve|^2)^{\frac{2-N}{2N-2}}\nabla \ve), &x\in\Omega, t\in(0,\Te)\\
    0 = \Delta \ve + \ue, &x\in\Omega, t\in(0,\Te)\\
    \partial_\nu \ue - \ue(\varepsilon+|\nabla \ve|^2)^{\frac{2-N}{2N-2}}\partial_\nu \ve = \ve = 0, &x\in\partial\Omega, t\in(0,\Te)\\
    \ue(\cdot,0) = u_0(\cdot), &x\in\Omega.
  \end{cases}
\end{equation}
Then the accumulated density $\Ue$ of $\ue$ solving  
\begin{equation*}
  \begin{cases}
    \Ue_t = N^2\xi^{2-\frac{2}{N}}\Ue_{\xi\xi} 
    + N(\varepsilon + \xi^{\frac2N-2} {\Ue}^2)^{\frac{2-N}{2N-2}}\Ue\Ue_\xi, 
    &\xi\in(0,1), t\in(0,\Te),\\
    \Ue(0,t) = 0, \quad \Ue(1,t) = \frac{m}{\omega_N}, &t\in(0,\Te),\\
    \Ue(\cdot,0) = U_0, &x\in(0,1),
  \end{cases}
\end{equation*}
is a subsolution of \eqref{sys: partial mass with flux limitation},
and thus Lemma~\ref{le: comparison principle} entails that  
$\Ue \leq U$ for all $(\xi,t)\in[0,1]\times[0,T^\varepsilon)$ and $\varepsilon \in (0,1)$ 
with $T^\varepsilon := \min\{T_{\max},\Te\}$.
We obtain from \eqref{eq: grad v =} and \eqref{h: U/s has upper bounds} that  
\begin{equation}
    \label{eq: boundedness of grad v}
  \begin{split}
    |\nabla \ve| 
    &= |\ve_r| 
    = r^{1-N}\int_0^r\ue\rho^{N-1}\dd \rho 
    \leq \frac{\Ue(r^N,t)}{r^N}\\
    &\leq \sup_{(\xi,t)\in(0,1)\times(0,T_{\max})} \frac{U(\xi,t)}{\xi}
    := C^{N-1}
  \end{split}
\end{equation}
for all $r\in(0,1)$, $t\in(0,T^\varepsilon)$ and $\varepsilon\in(0,1)$.
  Then 
  multiplying the first equation in \eqref{sys: regularized system} by $p{\ue}^{p-1}$ and integrating by parts, 
  we %
  have that for each $p\geq2$,
  \begin{equation}
    \label{eq: Lp estimate}
    \begin{split}
      \frac{\dd}{\dd t}\int_\Omega {\ue}^{p} 
      &= - \frac{4(p-1)}{p}\int_\Omega |\nabla {\ue}^{\frac{p}{2}}|^2 \\
      &\quad + p(p-1)\int_\Omega {\ue}^{p-1}(\varepsilon+|\nabla \ve|^2)^{\frac{2-N}{2N-2}}\nabla \ve \cdot \nabla \ue\\
      &\leq -2\int_\Omega |\nabla {\ue}^{\frac{p}{2}}|^2 
      + 2Cp\int_\Omega {\ue}^{\frac{p}{2}}|\nabla {\ue}^{\frac{p}{2}}|\\
      &\leq - \int_\Omega |\nabla {\ue}^{\frac{p}{2}}|^2 
      + 16C^2p^2\int_\Omega {\ue}^p
    \end{split}
  \end{equation}
  for all $t\in(0,T^\varepsilon)$ and $\varepsilon\in(0,1)$.
  Using the Gagliardo-Nirenberg inequality
  \begin{equation*}
  \|\psi\|_{L^{2+\frac{4}{Np}}}^{2+\frac{4}{Np}}
  \leq C_{\mathrm{GN}}\|\nabla \psi\|_{L^2(\Omega)}^2 
  \|\psi\|_{L^{\frac{2}{p}}}^{\frac{4}{Np}} 
  + C_{\mathrm{GN}}\|\psi\|_{L^{\frac{2}{p}}}^{2+\frac{4}{Np}} 
  \quad \text{for all } \psi\in W^{1,2}(\Omega),
  \end{equation*}
  with some $C_{\mathrm{GN}} > 0$ (independent of $p\geq2$) and Young inequality,
  we may infer from \eqref{eq: Lp estimate} that 
  \begin{equation*}
    \begin{split}
      \frac{\dd}{\dd t}\int_\Omega {\ue}^{p} &+ \int_\Omega {\ue}^{p}
      \leq - \int_\Omega |\nabla {\ue}^{\frac{p}{2}}|^2 
      + 16(C+1)^2p^2\int_\Omega {\ue}^p\\
      &\leq - \int_\Omega |\nabla {\ue}^{\frac{p}{2}}|^2 
      + C_{\mathrm{GN}}^{-1}m^{-\frac{2}{N}}\int_\Omega {\ue}^{p+\frac{2}{N}} 
      + 16^{1+\frac{Np}{2}}(C+1)^{2+Np}p^{2+Np}C_{\mathrm{GN}}^{\frac{Np}{2}}m^p|\Omega|\\
      &\leq m^p 
      + 16^{1+\frac{Np}{2}}(C+1)^{2+Np}p^{2+Np}C_{\mathrm{GN}}^{\frac{Np}{2}}m^p|\Omega|
    \end{split}
  \end{equation*}
  for all $t\in(0,T^\varepsilon)$ and $\varepsilon\in(0,1)$.
  It follows from the Gronwall inequality that for all $p\geq2$, there exists $C(p) > 0$ independent of $\varepsilon\in(0,1)$ such that 
  \begin{equation}
    \label{eq: u Lp integrable}
    \sup_{t\in(0,T^\varepsilon)}\int_\Omega {\ue}^p < C(p) 
    \quad\text{for all } \varepsilon \in(0,1).
  \end{equation}
  Thanks to the Moser iteration performed in \cite[Lemma~A.1]{Tao2012},
  one may invoke \eqref{eq: boundedness of grad v} and \eqref{eq: u Lp integrable} to find a constant $K>0$ independent of $\varepsilon\in(0,1)$ such that  
  \begin{equation}
    \label{eq: ue is uniformly bounded}
  \|\ue\|_{L^\infty(\Omega)} \leq K\quad \text{for all } t\in(0,T^\varepsilon),
  \end{equation}  
  which implies $\Te\geq T_{\max} = T^\varepsilon$ for all $\varepsilon\in(0,1)$.
  It follows from \eqref{eq: Lp estimate} with $p=2$ that for each $T\in(0,T_{\max})$, 
  $\ue$ and $\ue_t$ are uniform-in-$\varepsilon$ bounded in $L^2((0,T); W^{1,2}(\Omega))\cap C^0([0,T]; L^2(\Omega))$ and $L^2((0,T); W^{-1,2})$, respectively.
  Moreover, for any $T\in(0,T_{\max})$ and $q > N$,
  $\ve$ is uniform-in-$\varepsilon$ bounded in $C^0([0,T];W^{2,q}(\Omega)\cap W^{1,q}_0(\Omega))$.
  Therefore, by standard Aubin-Lions arguments as done in \cite{Yan2020,Kohatsu2024}, 
  we may extract a subsequence if necessary such that $(\ue,\ve)$ approaches to the unique weak solution $(u,v)$ of \eqref{sys: my ks flux limitation},
  as $\varepsilon\downarrow0$.
  The desired boundedness is a consequence of \eqref{eq: ue is uniformly bounded}.
\end{proof}

We need the following auxiliary lemma~\cite{Mao2024} to establish an upper bound for the initial data.

\begin{lemma}
  \label{le: function with bound derivative can be bounded}
  Suppose the function $f\in C^0([0,1])\cap C^1((0,1))$ is such that
  $f(0)=0$ 
  and 
   \begin{equation*}
     0 \leq \sup_{\xi\in(0,1)} f'(\xi) < \infty.
   \end{equation*}
  Let 
  $\{f_a: [0,1]\mapsto[0,\infty) \mid f_a(0)=0\}_{a\in \mathcal{I}}$
  be a family of concave and increasing functions satisfying 
  \begin{equation*}
  \label{eq: property of a family concove functions}
    \sup_{a\in\mathcal{I}} f_a(\xi) > f(1)
    \quad\text{for all } \xi\in(0,1).
  \end{equation*}
  Then there exists \(\alpha\in\mathcal{I}\) such that
  \[
  f(\xi) < f_\alpha(\xi)
  \quad\text{for all } \xi\in(0,1).
  \]
\end{lemma}

Now we are in a position to show Proposition~\ref{prop: varepsilon regualrity}.

\begin{proof}
  [Proof of Proposition~\ref{prop: varepsilon regualrity}]
  If the accumulated density of a radially symmetric solution to \eqref{sys: my ks flux limitation} satisfies \eqref{eq: U(rho,t) < A},
  then there exists $\Lambda\in(0,A)$ such that 
  $U \leq U(\rho) < \Lambda$ for all $(\xi,t)\in[0,\rho]\times[0,T_{\max})$.
Thanks to~\eqref{eq: phi ell uparrow A} and $U_0\in C^1([0,1])$ as $u_0\in C^0(\overline{\Omega})$, 
we may fix $\ell\in(\Lambda, A)$ according to Lemma~\ref{le: function with bound derivative can be bounded} such that 
  \begin{equation*}
  U_0(\rho\xi)\leq \phi_\ell(\xi) 
  \quad\text{for all } \xi\in(0,1),
  \end{equation*} 
  where $\phi_\ell$ is defined in Corollary~\ref{coro: steady states}.
  Then Lemma~\ref{le: comparison principle} assures that 
  \begin{equation*}
    U(\rho\xi) \leq \phi_\ell(\xi) 
    \quad \text{for all } (\xi,t)\in[0,1]\times[0,T_{\max}),
  \end{equation*}
  where we have used the scaling invariant of \eqref{sys: partial mass with flux limitation}.
  Therefore, 
  \begin{equation*}
    \frac{U(\xi,t)}{\xi} 
    \leq \max\bigg\{\frac{\phi_{\ell\xi}(0)}{\rho}, \frac{U_0(1)}{\rho}\bigg\}
    \quad \text{for all } (\xi,t)\in(0,1)\times[0,T_{\max}).
  \end{equation*}
  Lemma~\ref{le: varepsilon regularity} is applicable to assert that \eqref{eq: U xi has upper bounds} holds, 
  which completes the proof.  
\end{proof}

\subsection{Stability}
Proposition~\ref{prop: varepsilon regualrity} has already established that any radially symmetric solution with subcritical mass exists globally and remains bounded.
We now aim to show that such a globally bounded solution converges to its steady state $X_\lambda$ defined in \eqref{sym: steady states of Cauchy problem} and 
uniquely determined by initial mass.

\begin{theorem}\label{thm: global boundedness}
  Let $N\geq2$ and $\Omega = B_1\subset\mathbb{R}^N$. 
  For any choice of initial function $u_0\in C^0(\overline{\Omega})$ complying with \eqref{h: initial conditions} and 
  \begin{equation*}
    m = \int_\Omega u_0\dd x < m_c,
  \end{equation*}
the solution $(u,v)$ of \eqref{sys: my ks flux limitation} is globally bounded and approaches to the steady state $X_\lambda$ in the topology of $L^\infty(\Omega)$,
where $\lambda > 0$ is uniquely determined by $\int_\Omega u_0\dd x = \int_\Omega X_\lambda\dd x$.
\end{theorem}

We shall establish a Lyapronov functional. 
We pay our attention to the cases $N\geq3$ and refer to \cite{Biler2006} for the case of $N=2$.

\begin{lemma}
  \label{le: dissipation inequality}
  Let $N\geq2$ and $\Omega=B_1\subset\mathbb{R}^N$.
  Define 
  \begin{equation}
    \label{sym: Psi ell}
    \Psi_\ell := \int_0^1\xi^{\frac{2}{N}-1}(2-\xi)|R_\ell|\dd\xi 
    \quad \text{for } t\in[0,\infty),
  \end{equation}
  where 
  \begin{equation*}
    R_\ell := U - \phi_\ell
    \quad \text{for } (\xi,t)\in[0,1]\times[0,\infty),
  \end{equation*}
  and $\phi_\ell$ is defined in Corollary~\ref{coro: steady states}
  with $\ell = U_0(1)\in(0,A)$.
  Then 
  \begin{equation}
    \label{eq: dissipation of R}
    \begin{split} 
    \Psi_\ell'
    &\leq  \frac{(2-N)N^2}{N-1}\int_0^1|R_\ell|\dd\xi
    \end{split}
  \end{equation}
  for all $t\in(0,\infty)$.
\end{lemma}

\begin{proof}
  The proof is a minor modification of the stability argument in \cite[Theorem 3.1]{Biler2006}.
  Recall \eqref{sym: Pd} and define 
\begin{equation}
  \Phi_\delta(r) := \begin{cases}
    \Pd(r) &\text{if } r\geq0,\\
    \Pd(-r) &\text{if } r < 0,
  \end{cases}
\end{equation}
which is a convex approximation of $r\mapsto|r|$.
Indeed, $r\mapsto\Phi_\delta(r)$ and $r\mapsto r\Phi_\delta'(r)$ converge uniformly to $|r|$ over $\mathbb{R}$, and $r\mapsto r\Phi_\delta''(r)$ is bounded and converges a.e. to $0$ as $\delta\downarrow0$.

Let $\varrho := 2 - \xi$.
Multiplying \eqref{sym: mathcal P} by $\varrho\Phi_\delta'(R_\ell)\xi^{\frac{2}{N}-1}$ and integrating over $(0,1)$,
we obtain 
\begin{align} 
  \frac{\dd}{\dd t}\int_0^1\varrho\Phi_\delta(R_\ell)\xi^{\frac{2}{N}-1}\dd\xi 
  &= \int_0^1\bigg(N^2\xi U_\xi - N^2U + (N-1)U^{\frac{N}{N-1}}\bigg)_\xi\varrho\Phi_\delta'(R_\ell)\dd\xi\notag\\
  &=  \int_0^1\bigg(N^2\xi U_\xi - N^2U + (N-1)U^{\frac{N}{N-1}}\bigg)\Phi_\delta'(R_\ell)\dd\xi\notag\\
  &\quad - \int_0^1\bigg(N^2\xi U_\xi - N^2U + (N-1)U^{\frac{N}{N-1}}\bigg)\varrho\Phi_\delta''(R_\ell)R_{\ell\xi}\dd\xi\notag\\
  &=: I_1 - I_2 \label{eq: d/dt Phi delta R ell}
\end{align} 
for all $t\in(0,\infty)$, 
because of $\Phi_\delta'(0) = 0$ and $U(0,t) = R_\ell(1,t) = 0$ for all $t>0$.
Using the identity \eqref{sys: ode of W}, i.e., 
\begin{equation*}
  0 = N^2\xi\phi_{\ell\xi\xi} + N\phi_\ell^{\frac{1}{N-1}}\phi_{\ell\xi} 
  = (N^2\xi\phi_{\ell\xi} - N^2\phi_\ell + (N-1)\phi_\ell^{\frac{N}{N-1}})_\xi
  \quad \text{for } \xi\in(0,1),
\end{equation*}
we get by integrating by parts
\begin{align*}
  I_1
  &= \int_0^1\bigg(N^2\xi R_{\ell\xi} - N^2R_\ell + (N-1)U^{\frac{N}{N-1}} - (N-1)\phi_\ell^{\frac{N}{N-1}}\bigg)\Phi_\delta'(R_\ell)\dd\xi\\
  &= - 2N^2\int_0^1R_\ell \Phi'_\delta(R_\ell)\dd\xi 
  - N^2\int_0^1\xi R_\ell\Phi''_\delta(R_\ell)R_{\ell\xi}\dd\xi \\
  &\quad + (N-1)\int_0^1 \bigg(U^{\frac{N}{N-1}} - \phi_\ell^{\frac{N}{N-1}}\bigg)\Phi_\delta'(R_\ell)\dd\xi
\end{align*}
for all $t\in(0,\infty)$, 
and 
\begin{align*}
  -I_2 &= - \int_0^1\bigg(N^2\xi R_{\ell\xi} - N^2R_\ell 
  + (N-1)U^{\frac{N}{N-1}} - (N-1)\phi_\ell^{\frac{N}{N-1}}\bigg)\varrho\Phi_\delta''(R_\ell)R_{\ell\xi}\dd\xi\\
  &\leq N^2\int_0^1\varrho R_\ell\Phi''_\delta(R_\ell)R_{\ell\xi}\dd \xi 
  - (N-1)\int_0^1\bigg(U^{\frac{N}{N-1}}- \phi_\ell^{\frac{N}{N-1}}\bigg)\varrho\Phi_\delta''(R_\ell)R_{\ell\xi}\dd\xi
\end{align*}
for all $t\in(0,\infty)$ because of $\Phi_\delta''\geq0$ for a.e. $r\in\mathbb{R}$.
Therefore, \eqref{eq: d/dt Phi delta R ell} is reduced to 
\begin{align*}
  \frac{\dd}{\dd t}\int_0^1\varrho\Phi_\delta(R_\ell)\xi^{\frac{2}{N}-1}\dd\xi
  &\leq - 2N^2\int_0^1R_\ell \Phi_\delta'(R_\ell)\dd\xi 
  + (N-1)\int_0^1 \bigg(U^{\frac{N}{N-1}} - \phi_\ell^{\frac{N}{N-1}}\bigg)\Phi_\delta'(R_\ell)\dd\xi \\
  &\quad + N^2\int_0^1\varrho R_\ell\Phi''_\delta(R_\ell)R_{\ell\xi}\dd \xi 
      - N^2\int_0^1\xi R_\ell\Phi_\delta''(R_\ell)R_{\ell\xi}\dd\xi \\
  &\quad - (N-1)\int_0^1\bigg(U^{\frac{N}{N-1}}- \phi_\ell^{\frac{N}{N-1}}\bigg)\varrho\Phi_\delta''(R_\ell)R_{\ell\xi}\dd\xi
\end{align*}
for all $t\in(0,\infty)$. 
Thus, Lebesgue dominated convergence theorem ensures we can pass to the limit as $\delta \downarrow 0$ in the above inequality, and end up with 
\begin{equation}
  \label{eq: dissipation inequality}
  \begin{split} 
  \Psi_\ell' 
  &\leq - 2N^2\int_0^1 |R_\ell|\dd\xi   
  + (N-1)\int_0^1\bigg\lvert U^{\frac{N}{N-1}} - \phi_\ell^{\frac{N}{N-1}}\bigg\rvert\dd\xi
  \end{split}
\end{equation}
for all $t\in(0,T_{\max})$.
We estimate the last term of \eqref{eq: dissipation inequality} 
by virtue of %
mean value theorem~\eqref{eq: mean value theorem}, 
\begin{equation*}
  (N-1)\int_0^1\bigg\lvert U^{\frac{N}{N-1}} - \phi_\ell^{\frac{N}{N-1}}\bigg\rvert\dd\xi 
  \leq NA^{\frac{1}{N-1}}\int_0^1|R_\ell|\dd\xi
  = \frac{N^3}{N-1}\int_0^1|R_\ell|\dd\xi,
\end{equation*}
which implies \eqref{eq: dissipation of R} holds.
\end{proof}

\begin{remark}
  In the case of $N=2$, we have $\ell < A=4$. 
  So it is easy to see from \eqref{eq: dissipation inequality} that 
  \begin{equation*}
    \Psi'_\ell \leq - 8\int_0^1 |R_\ell|\dd\xi   
    + \int_0^1 (U + \phi_\ell)\lvert R_\ell\rvert\dd\xi
    \leq - (8-2\ell)\int_0^1|R_\ell|\dd\xi 
    \quad \text{for all } t\in(0,\infty),
  \end{equation*}
  which has been shown in~\cite{Biler2006} and proven to be sufficient for \eqref{eq: Psi ell to 0} below.
\end{remark}

Now we are in a position to show Theorem~\ref{thm: global boundedness}.

\begin{proof}
  [Proof of \upshape{Theorem~\ref{thm: global boundedness}}]
Global boundedness has been established in Proposition~\ref{prop: varepsilon regualrity}. 
Now we show globally asymptotic stability and ignore the case of $N=2$ that has been shown in~\cite{Biler2006}.
Lemma~\ref{le: dissipation inequality} warrants the function $\Psi_\ell$ is nonincreasing.
  We claim 
\begin{equation}
  \label{eq: Psi ell to 0}
  \Psi_\ell\downarrow0,
  \quad \text{as } t\to\infty.
\end{equation}
Otherwise, there exists $C>0$ such that 
\begin{equation*}
  \Psi_\ell > C 
  \quad\text{for all } t\in(0,\infty).
\end{equation*}
As 
\begin{equation*}
  \int_0^\eta\xi^{\frac{2}{N}-1}(2-\xi)|R_\ell|\dd\xi 
  < NA\eta^{\frac{2}{N}}
  \quad\text{for all } (\eta,t)\in(0,1)\times(0,\infty),
\end{equation*}
we may fixed $\eta\in(0,1)$ small enough such that  
\begin{equation*}
  \int_\eta^1\xi^{\frac{2}{N}-1}(2-\xi)|R_\ell|\dd\xi 
  \geq \Psi_\ell - NA\eta^{\frac{2}{N}} 
  > \frac{C}{2} 
  \quad\text{for all } t\in(0,\infty),
\end{equation*}
which implies 
\begin{equation*}
  \int_\eta^1|R_\ell|\dd\xi > \frac{C}{4}\eta^{1-\frac{2}{N}} 
  \quad\text{for all } t\in(0,\infty).
\end{equation*}
It follows from \eqref{eq: dissipation of R} that 
\begin{equation*}
  \begin{split} 
  \Psi_\ell &\leq \Psi_\ell(0) - \frac{(N-2)N^2}{N-1}\int_0^t\int_0^1|R_\ell|\dd\xi\dd\tau\\
  &\leq NA - \frac{(N-2)N^2C\eta^{1-\frac{2}{N}}}{4(N-1)}t
  \quad \text{for all } t\in (0,\infty),
  \end{split}
\end{equation*}
and thus $\Psi_\ell$ becomes negative in finite time. 
This contradicts the nonnegativity of $\Psi_\ell$ according to its definition~\eqref{sym: Psi ell}, as desired.

Thanks to \eqref{eq: Psi ell to 0}, we immediately obtain  
\begin{equation*}
  \int_0^1|R_\ell|\dd\xi \to 0,
  \quad \text{as } t\to\infty.
\end{equation*}
Using 
\begin{equation*}
  \|u\|_{C^{\vartheta,\frac{\vartheta}{2}}(\overline{\Omega}\times[1,\infty))} < \infty,
\end{equation*} 
for some $\vartheta\in(0,1)$,
according to \cite[Thoerem~1.3]{Porzio1993},
we have 
\begin{equation*}
  \frac{|U_\xi(\xi,t)-U_\eta(\eta,t)|}{|\xi-\eta|^{\frac{\vartheta}{N}}} 
  \leq \frac{|u(\xi^{\frac{1}{N}},t)-u(\eta^{\frac{1}{N}},t)|}{N|\xi-\eta|^{\frac{\vartheta}{N}}}
  \leq \frac{|u(\xi^{\frac{1}{N}},t)-u(\eta^{\frac{1}{N}},t)|}{N|\xi^{\frac{1}{N}}-\eta^{\frac{1}{N}}|^\vartheta} 
  \leq \|u\|_{C^{\vartheta,\frac{\vartheta}{2}}(\overline{\Omega}\times[1,\infty))}
\end{equation*}
for all $\xi,\eta\in(0,1)$, $t>1$ and $\xi\neq\eta$.
So $\|U\|_{C^{1+\frac{\vartheta}{N}}([0,1])} \leq  \|u\|_{C^{\vartheta,\frac{\vartheta}{2}}(\overline{\Omega}\times[1,\infty))}$ for all $t>1$. 
By an Ehrling type argument, 
we have that
\begin{equation*}
  R_\ell\to 0
  \quad\text{in } W^{1,\infty}((0,1)), 
  \text{ as } t\to\infty,
\end{equation*}
because of $C^{1+\frac{\vartheta}{N}}((0,1))\hookrightarrow \hookrightarrow W^{1,\infty}((0,1))\hookrightarrow L^1((0,1))$.
This, for $\lambda$ determined by \eqref{eq: lambda determined by ell}, implies 
\begin{equation*}
  u\to X_\lambda\quad\text{in } L^\infty(\Omega),
  \text{ as } t\to\infty.
\end{equation*}
The proof is complete.
\end{proof}

\section{Infinite-time chemotactic collapse}
\label{sec: infinite-time blow-up}

This section is devoted to show that for arbitrary radially symmetric initial datum with critical mass, 
the solution exists globally and exhibit a chemotactic collapse in infinite time.

\begin{theorem}
  \label{thm: infinite-time chemotactic collapse}
  Let $N\geq2$ and $\Omega = B_1\subset\mathbb{R}^N$.
  For any choice of initial function $u_0\in C^0(\overline{\Omega})$ complying with \eqref{h: initial conditions} and 
  \begin{equation*}
    \int_\Omega u_0\dd x = m_c,
  \end{equation*}
the corresponding solution $(u,v)$ of \eqref{sys: my ks flux limitation} exists globally and aggregates complete mass at the center of the domain in infinite time.
\end{theorem}

We adapt the techniques from~\cite{Nagai2000} to show that if radially symmetric initial data explodes in finite time, 
then the solution concentrates at least the critical mass $m_c$ at the origin.

\begin{lemma}
  \label{le: chemotactic collapse}
  Suppose a radially symmetric solution $(u,v)$ explodes in finite time $T_{\max}\in(0,\infty)$.
Then there exists a nonnegative and increasing function $U^\star\in C^2((0,1])$ such that  
\begin{equation*}
  U \to U^\star \quad\text{in } C^2_{\loc}((0,1]), \quad\text{as } t\uparrow T_{\max},
\end{equation*}
but 
\begin{equation}
  \label{eq: U(0+) geq A}
  \lim_{\xi\downarrow0}U^\star \geq A.
\end{equation}
\end{lemma}

\begin{proof}
  Using Helly's compact theorem, we may find a sequence $\{t_k\}$ and an increasing function $U^\star$ such that $t_k\uparrow T_{\max}$ as $k\to\infty$ and  
  \begin{equation*}
    U(\cdot, t_k) \to U^\star(\cdot) 
    \quad\text{for all } \xi\in(0,1], \quad\text{as } k\to\infty.
  \end{equation*}
  For $\delta\in(0,1)$, we may apply standard parabolic regularity theories to deduce that 
  \begin{equation*}
    \|U\|_{C^{2+\vartheta,1+\frac{\vartheta}{2}}((\delta,1]\times(T_{\max}/2,T_{\max}))} 
    < C(\delta,\vartheta),
  \end{equation*}
  for some positive constant $C(\delta,\vartheta) > 0$ and $\vartheta\in(0,1)$.
  So 
  \begin{equation*}
    \lim_{t\uparrow T_{\max}} U(\cdot, t) \text{ exists and equals to } U^\star(\cdot) 
    \quad\text{for all } \xi\in(0,1], %
  \end{equation*}
  which also holds in the sense of $C^2([\delta ,1])$ for all $\delta \in(0,1)$.
  If \eqref{eq: U(0+) geq A} is violated, 
  there exists $\varepsilon > 0$ such that 
  \begin{equation*}
    U^\star(\varepsilon) < A,
  \end{equation*}
  and thus 
  \begin{equation*}
    \sup_{t\in(T,T_{\max})} U(\varepsilon,t) < A
  \end{equation*}
  for some $T\in[0,T_{\max})$.
  Proposition~\ref{prop: varepsilon regualrity} entails 
  \begin{equation*}
    \|U_s\|_{L^\infty((0,1)\times(T,T_{\max}))} < \infty,
  \end{equation*}
  which is incompatible with the extensibility criterion~\eqref{eq: extensibility criterion}.%
\end{proof}

Following Suzuki~\cite{Suzuki2013}, 
we show that if radially symmetric initial data has critical mass, 
then the solution exists globally.

\begin{lemma}
  \label{le: global existence}
  Let $N\geq2$ and $\Omega=B_1\subset\mathbb{R}^N$, 
  and assume that $u_0\in C^0(\overline{\Omega})$ is nonnegative and radially symmetric. 
  If $\int_\Omega u_0\dd x = m_c$, 
  then the solution $(u,v)$ exists globally.
\end{lemma}

\begin{proof}
  Suppose by contradiction that the solution $(u,v)$ blows up in finite time $T_{\max}\in(0,\infty)$.
  Then it follows from Lemma~\ref{le: chemotactic collapse} that the solution $U$ of \eqref{sys: partial mass with flux limitation} can be continuously extended to $(0,1]\times\{T_{\max}\}$ and 
  \begin{equation*}
    U^\star(\xi) \equiv A \quad \text{for all } \xi\in(0,1].
  \end{equation*}
  This implies 
  \begin{equation*}
    U \equiv A \quad \text{for all } (\xi,t)\in(0,1]\times[0,T_{\max}),
  \end{equation*} 
  by strong maximal principle,
  which is incompatible with our assumption $u_0\in C^0(\overline{\Omega})$.
\end{proof}

We shall derive a dissipation functional for initial data with critical mass, 
analogous to Lemma~\ref{le: dissipation inequality}.

\begin{lemma}
  \label{le: dissipation inequality for critical case}
  Let $N\geq2$ and $\Omega=B_1\subset\mathbb{R}^N$.
  Assume that the initial datum $u_0\in C^0(\overline{\Omega})$ is a radially symmetric function with $\int_\Omega u_0\dd x = m_c$. 
  Define 
  \begin{equation*}
    \Psi_c := \int_0^1\xi^{\frac{2}{N}-1}(2-\xi)R_c\dd\xi 
    \quad \text{for } t\in(0,\infty),
  \end{equation*}
  where 
  \begin{equation*}
    R_c := A - U 
    \quad \text{for } (\xi,t)\in[0,1]\times[0,\infty).
  \end{equation*}
  Then 
  \begin{equation}
    \label{eq: dissipation of R_c}
    \begin{split} 
    \Psi_c'
    \leq -  \frac{(N-2)N^2}{N-1}\int_0^1R_c\dd\xi
    \end{split}
  \end{equation}
  for all $t\in(1,\infty)$.
\end{lemma}

\begin{proof}
  We proceed as in the proof of Lemma~\ref{le: dissipation inequality} and \cite[Proposition~3.2]{Biler2006}.
  We notice that $R_c$ solves
  \begin{equation}
    \label{eq: R ct}
    R_{ct} = - N^2\xi^{2-\frac{2}{N}}U_{\xi\xi} - N\xi^{1-\frac{2}{N}}U^{\frac{1}{N-1}}U_\xi
    = - \xi^{1-\frac{2}{N}}\bigg(N^2\xi U_\xi - N^2U + (N-1)U^{\frac{N}{N-1}}\bigg)_\xi.
  \end{equation}
 Multiplying \eqref{eq: R ct} by $\xi^{\frac{2}{N}-1}(2-\xi)$ and integrating over $(0,1)$,
 we obtain 
\begin{align}
  \Psi'_c &= - \int_0^1(2-\xi)\bigg(N^2\xi U_\xi - N^2 U + (N-1)U^{\frac{N}{N-1}}\bigg)_\xi\dd\xi\notag\\
  &= -N^2U_\xi(1,t) - \int_0^1\Big(N^2\xi U_\xi -N^2U + (N-1)U^{\frac{N}{N-1}}\Big)\dd\xi\notag\\
  &= - N^2U_\xi(1,t) 
  - N^2A 
  + 2N^2\int_0^1U\dd\xi 
  - (N-1)\int_0^1U^{\frac{N}{N-1}}\dd\xi
  \label{eq: Psi'c}
\end{align}
because of $U(1,t) = A$ and $U(0,t) = 0$ for all $t\in(0,\infty)$.
Noting  
\begin{equation*}
  2N^2\int_0^1U\dd\xi = -2N^2\int_0^1R_c\dd\xi + 2N^2A,
\end{equation*}
and in view of mean value theorem~\eqref{eq: mean value theorem}
\begin{align*}
  -(N-1)\int_0^1U^{\frac{N}{N-1}}\dd\xi 
  &= (N-1)\int_0^1\bigg(A^{\frac{N}{N-1}}-U^{\frac{N}{N-1}}\bigg)\dd\xi - (N-1)A^{\frac{N}{N-1}}\\
  &\leq NA^{\frac{1}{N-1}}\int_0^1R_c\dd\xi - (N-1)A^{\frac{N}{N-1}}\\
  &= \frac{N^3}{N-1}\int_0^1R_c\dd\xi - N^2A,
\end{align*}
we obtain \eqref{eq: dissipation of R_c}.
\end{proof}

\begin{remark}
  In the case of $N=2$, we have $A=4$. 
  So it is easy to see from \eqref{eq: Psi'c} that 
  \begin{equation*}
    \Psi'_c \leq - \int_0^1R_c^2\dd\xi 
    \leq - \frac{1}{4}\int_0^1(2-\xi)^2R_c^2\dd\xi
    \leq - \frac{1}{4}\Psi_c^2
    \quad \text{for all } t\in(0,\infty),
  \end{equation*}
  by H\"older inequality, which has been shown in~\cite{Biler2006} and proven to be sufficient for \eqref{eq: Rc tends to 0}.
\end{remark}

Now we are in a position to show Theorem~\ref{thm: infinite-time chemotactic collapse}.

\begin{proof}
  [Proof of \upshape{Theorem~\ref{thm: infinite-time chemotactic collapse}}]
  Following the lines in the proof of Theorem~\ref{thm: global boundedness},
  we have 
  \begin{equation}
    \label{eq: Rc tends to 0}
    \int_0^1R_c\dd\xi \to 0,
    \quad\text{as } t\to\infty,
  \end{equation}
  according to Lemma~\ref{le: dissipation inequality for critical case}.
  We claim 
  \begin{equation}
    \label{eq: u to mc delta}
    u(\cdot,t)\weakstar m_c\delta_0(\cdot) 
    \quad\text{in } \mathcal{M}(\overline{\Omega})
    \text{ as } t\to\infty.
  \end{equation}
  Suppose by absurdum that \eqref{eq: u to mc delta} is violated along a unbounded sequence $\{t_k\}_{k\in\mathbb{N}}$.
  Then by Banach-Alaoglu theorem, 
  we may extract a subsequence if necessary such that  
  \begin{equation}
    \label{eq: mu neq mcdelta}
    u(\cdot, t_k)\weakstar \mu(\cdot) \neq m_c\delta_0 
    \quad\text{in } \mathcal{M}(\overline{\Omega})
  \end{equation}
  for some positive measure $\mu\in\mathcal{M}(\overline{\Omega})$.
  Let $\eta_\delta(r) := \eta(\delta^{-\frac{1}{N}}r)$ for $\delta > 0$, 
  where 
  \begin{equation*}
    \eta (r) := 
    \begin{cases}
      1 & \text{if } r\in(-\infty,1),\\
      2-r & \text{if } r\in[1,2],\\
      0 & \text{if } r\in(2,\infty).
    \end{cases}
  \end{equation*}
  Define 
  \begin{equation*}
    f_k(\xi) := \int_{\Omega} u(\cdot,t_k)\eta_\xi(|x|)\dd x 
    \quad \text{and} \quad 
    f(\xi) := \int_\Omega \eta_\xi(|x|)\dd\mu(x) 
    \quad \text{for } \xi\in(0,1).
  \end{equation*}
  Then increasing and continuous function $f_k$ converges to increasing function $f$ pointwise, as $k\to\infty$,
  which implies 
  \begin{equation}
    \label{eq: fk to f}
    f_k \to f \quad\text{in } L^1((0,1)) \quad \text{ as } k\to\infty,
  \end{equation}
  by dominated convergence theorem.
  Since 
  \begin{equation*}
    \omega_NU(\xi,t_k) \leq f_k(\xi) \leq m_c
    \quad\text{for all } \xi\in(0,1) \text{ and } k\in\mathbb{N},
  \end{equation*}
  we get 
  \begin{equation*}
    \int_0^1|f_k-m_c|\dd\xi \leq \omega_N\int_0^1R_c(\cdot, t_k)\dd\xi
    \quad \text{for all } k\in\mathbb{N}.
  \end{equation*}
  It follows from \eqref{eq: Rc tends to 0} that  
  \begin{equation}
    \label{eq: fk to mc}
    f_k \to m_c \quad\text{in } L^1((0,1)) \quad \text{ as } k\to\infty.
  \end{equation}
  Hence, in view of \eqref{eq: fk to f} and \eqref{eq: fk to mc}, 
  we have $f \equiv m_c$ for $\xi\in(0,1)$ by monotonicity.
  We estimate 
  \begin{equation*}
    \begin{split}
      &\quad \bigg|\int_\Omega u(\cdot,t_k)\phi\dd x - m_c\phi(0)\bigg| \\
      &\leq \bigg|\int_\Omega u(\cdot,t_k)\eta_\delta(|x|)\phi(0)\dd x - m_c\phi(0)\bigg|
      + \bigg|\int_\Omega u(\cdot,t_k)\eta_\delta(|x|)(\phi-\phi(0))\dd x\bigg|\\
      &\quad + \bigg|\int_\Omega u(\cdot,t_k)(1-\eta_\delta(|x|))\phi\dd x\bigg|\\
      &\leq 2\bigg|\int_\Omega u(\cdot,t_k)\eta_\delta(|x|)\dd x - m_c\bigg|\|\phi\|_{L^\infty(\Omega)} 
      + m_c\|\eta_\delta(|x|)(\phi-\phi(0))\|_{L^\infty(\Omega)}
    \end{split}
  \end{equation*}
  for all $k\in\mathbb{N}$, $\delta > 0$ and $\phi\in C^0(\overline{\Omega})$.
  Therefore, by pointwise convergence of $f_k$ to $m_c$,
  \begin{equation*}
    \limsup_{k\to\infty} 
    \bigg|\int_\Omega u(\cdot,t_k)\phi\dd x - m_c\phi(0)\bigg| 
    \leq m_c\|\eta_\delta(|x|)(\phi-\phi(0))\|_{L^\infty(\Omega)}
  \end{equation*}
  for all $\delta > 0$ and $\phi\in C^0(\overline{\Omega})$. 
  Pushing $\delta\downarrow0$, 
  we have 
  \begin{equation*}
    \int_\Omega u(\cdot,t_k)\phi\dd x \to m_c\phi(0) 
    \quad \text{as } k\to\infty
  \end{equation*}
  for all $\phi\in C^0(\overline{\Omega})$,
  which is incompatible with our assumption~\eqref{eq: mu neq mcdelta}.
\end{proof}

\section*{Acknowledgment}

The first author is supported by the Fundamental Research Funds for the Central Universities (No.~B250201215). The third author is supported by  Natural Science Foundation of Jiangsu Province (No. SBK2022044966).

\end{document}